\newcommand{\qed}{\vbox{\hrule width4pt height3pt depth1pt} }

\newcommand{\emp}[1]{\bf #1\rm}

\font\bbb=msbm10
\newcommand{\R}{\hbox{\bbb\char82}}

\newcommand{\del}{\partial}

\newtheorem{theorem}{Theorem}

\newenvironment{proof}{{\bf Proof: }}{\hfill$\qed$}

\newcommand{\V}{{\bf V}}
\newcommand{\F}{{\bf F}}

\newcommand{\x}{{\bf x}}

\newcommand{\vz}{{\bf 0}}

\newcommand{\cF}{{\cal F}}

\newcommand{\id}{{\bf I}}

\documentclass{article}

\usepackage{graphicx}

\title{Cyclic Approximation to K-Stasis}
\author{Stewart D. Johnson\\ \\Department of Mathematics and Statistics\\Williams College}

\renewcommand{\arraystretch}{1.5}

\begin{document}             
\openup 4pt
\maketitle
\begin{abstract}
\noindent 
\end{abstract}

If a linear combination of k smooth vector fields is zero at a point, then, generically, near this point there are small cycles comprised of segments from the flow of each vector field. This answers a question posed in arXiv:math/0504365.  

\bigskip


\vfil\break
%

\section{Definitions and Results}

Vector fields $$\V_j:\R^n\to\R^n \hbox{ for } j=1,\ldots,k$$ induce flows $$\F_j(\x,t):\R^n\times \R \to \R^n.$$ 

A point $\x_0\in \R^n$ is a \emp{K-stasis point} if $$\sum_{j=1}^k m_j \V_j(\x_0)=\vz$$ for some probability weighting $(m_1,\ldots, m_k)$; where each $m_j > 0$ and $\sum m_j=1$. 

The stasis point is \emp{regular} if $$ \sum_{j=1}^k  m_j {\del \V_j \over \del \x} (\x_0)$$ is non-singular.  

A \emp{K-cycle} is a sequence of points $\x_1,\ldots,\x_k$ in $\R^n$, and a sequence of times $(\delta_1,\ldots,\delta_k)$ with each $\delta_j>0$, such that 
\begin{eqnarray*} 
\F_1(\x_1,\delta_1) &=& \x_2 \cr  
\F_2(\x_2,\delta_2) &=& \x_3 \cr
&\vdots& \cr
\F_{k}(\x_{k},\delta_{k}) &=& \x_1.
\end{eqnarray*}

We have the following theorem.

\begin{theorem} 
If $\x_0$ is a regular K-stasis point with weighting $(m_1,\ldots,m_k)$ then for all sufficiently small $\delta>0$ there exists a K-cycle with the time vector $(\delta m_1, \ldots, \delta m_k)$. 
 
\end{theorem}

\vfill\break

\begin{proof}

Without loss of generality take $\x_0=0$. Define
$$ \cF:\underbrace{\R^n\times\cdots\times\R^n}_{k} \times \R \to \R^n$$
as an average velocity
$$\cF(\x_1,\ldots,\x_k,\delta) = {(\F_1(\x_1,\delta m_1)-\x_1)\; +  \cdots  + \; (\F_k(\x_k,\delta m_k)-\x_k) \over  \delta}.$$

Note that $\cF$ is $C^1$ near $(\vz,\ldots,\vz,0)$, and that for $\delta=0$,  $${\del\over\del \x_j}\cF(\x_1,\ldots,\x_k,0) = m_j {\del \V_j \over \del \x} (\x_j)$$

Now $$\left(
\begin{array}{c}
\cF(\x_1,\ldots,\x_k,\delta)\cr 
\F_1(\x_1,\delta m_1) - \x_2 \cr 
\F_2(\x_2,\delta m_2) - \x_3 \cr 
\vdots \cr 
\F_{k-1}(\x_{k-1},\delta m_{k-1}) - \x_k
\end{array}
\right): \underbrace{\R^n\times\cdots\times\R^n}_{k} \times \R \to \underbrace{\R^n\times\cdots\times\R^n}_{k}$$
with
$$\left.\left(
\begin{array}{c}
\cF(\x_1,\ldots,\x_k,\delta)\cr 
\F_1(\x_1,\delta m_1) - \x_2 \cr 
\vdots \cr 
\F_{k-1}(\x_{k-1},\delta m_{k-1}) - \x_k
\end{array}
\right)\right|_{(\vz,\ldots,\vz,0)} 
= \left(
\begin{array}{c}
\vz\cr\vdots\cr\vz
\end{array}\right)$$

By the implicit function theorem, 
$$\left.\left(
\begin{array}{c}
\cF(\x_1,\ldots,\x_k,\delta)\cr 
\F_1(\x_1,\delta m_1) - \x_2 \cr 
\vdots \cr \F_{k-1}(\x_{k-1},\delta m_{k-1}) - \x_k
\end{array} 
\right)\right|_{(\x_1,\ldots,\x_k,\delta)} 
= \left(
\begin{array}{c}
\vz\cr\vdots\cr\vz
\end{array}\right)$$
will have solutions 
$$\x_1(\delta),\ldots, \x_k(\delta)$$ for small non-zero $\delta$ provided that the $nk\times nk$ matrix 
\renewcommand{\arraystretch}{2}
$$\left[
\begin{array}{c}
{\del  \cF (\x_1,\ldots,\x_k,\delta)\over \del \x_1,\ldots,\x_k} \cr
{\del (\F_1(\x_1,\delta m_1) - \x_2)\over \del \x_1,\ldots,\x_k}  \cr 
\vdots \cr 
{\del (\F_{k-1}(\x_{k-1},\delta m_{k-1}) - \x_k)\over \del \x_1,\ldots,\x_k}  
\end{array}
\right]_{(\vz,\ldots,\vz,0)}$$ is non-singular. 

\renewcommand{\arraystretch}{1.5}
This evaluates to
$$\left[ 
\begin{array}{ccccc} 
m_1{\del \V_1 \over \del \x} & m_2{\del \V_2 \over \del \x} & m_3{\del \V_3 \over \del \x} & \cdots & m_k{\del \V_k \over \del \x}\cr
\id & -\id &  \vz & \cdots & \vz \cr
\vz &  \id & -\id & \cdots & \vz \cr
\vdots & \vdots & \vdots & \ddots & \vdots \cr
\vz &  \vz &  \vz & \cdots & -\id \cr
\end{array}
\right]_{(\vz,\ldots,\vz,0)}
$$
This matrix is singular iff $\sum_{j=1}^k  m_j {\del \V_j \over \del \x}(\vz)$ is singular, and the result follows.  
\end{proof}

%
%
%
%
%
\end{document}